\documentclass{amsart}
\usepackage{graphics,url}
\renewcommand{\Re}{\mathbb{R}}
\newcommand{\Int}{\operatorname{Int}}
\newtheorem{thm}{Theorem}[section]
\begin{document}
{\noindent\small 
Topology Atlas Invited Contributions {\bf 9} (2004) no.~1, 12 pp.}
\bigskip

\title[Computational topology for regular closed sets]{Computational 
topology for regular closed sets (within the I-TANGO Project)}
\author[I-TANGO Project]{T.J.~Peters}
\address{T.J.~Peters, University of Connecticut, Storrs, CT 06269-3155 USA}
\email{tpeters@cse.uconn.edu.}
\author[]{J.~Bisceglio}
\address{J.~Bisceglio, Bluesky Studios, Inc.}
\author[]{D.R.~Ferguson}
\address{D.R.~Ferguson, The Boeing Company, Seattle, Washington}
\email{david.r.ferguson@boeing.com}
\author[]{C.M.~Hoffmann}
\address{C.M.~Hoffmann, Purdue University, East Lafayette, IN}
\email{cmh@cs.purdue.edu}
\author[]{T.~Maekawa}
\address{T.~Maekawa, Yokohama National University, Yokohama 240-8501 
Japan \& MIT}
\email{maekawa@ynu.ac.jp}
\author[]{N.M.~Patrikalakis}
\address{N.M.~Patrikalakis, Massachusetts Institute of Technology, 
Cambridge, MA}
\email{nmp@mit.edu}
\author[]{T.~Sakkalis}
\address{T.~Sakkalis, Agricultural University of Athens, Athens, Greece 
\& MIT}
\email{stp@aua.gr}
\author[]{N.F.~Stewart}
\address{N.F.~Stewart, Universit\'e de Montr\'eal, Montr\'eal, Canada}
\email{stewart@iro.umontreal.ca}
\thanks{Partial funding for all authors was from NSF grant DMS-0138098.  
Authors Peters and Sakkalis were also partially supported by NSF grant CCR 
0226504.  
All statements in this publication are the responsibility of the authors, 
\emph{not} of any of these indicated funding sources.}
\date{November 2003}

\begin{abstract}
The Boolean algebra of regular closed sets is prominent in topology, 
particularly as a dual for the Stone-$\check{\mathrm{C}}$ech 
compactification.  
This algebra is also central for the theory of geometric computation, as 
a representation for combinatorial operations on geometric sets.  
However, the issue of computational approximation introduces unresolved 
subtleties that do not occur within ``pure'' topology.  
One major effort towards reconciling this mathematical theory with 
computational practice is our ongoing I-TANGO project.  
The acronym I-TANGO is an abbreviation for ``Intersections---Topology, 
Accuracy and Numerics for Geometric Objects''.  The long-range goals and 
initial progress of the I-TANGO team in development of computational topology 
are presented.
\end{abstract}
\maketitle

\section{Introduction and brief literature review}

Throughout this paper, all sets considered will be assumed to be subsets
of $\Re^3$, with its usual topology.  The Boolean algebra of regular
closed sets in $\Re^3$ will be denoted as $\mathcal{R}(\Re^3)$.  
Furthermore, any regular closed set considered will be assumed to be
compact.  Any surfaces and curves considered will be assumed to be compact
2-manifolds and 1-manifolds, respectively.  All neighborhoods will be
assumed to be open subsets of $\Re^3$.

\begin{figure}[h]%[htbp]
%\centerline{\rotatebox{0}{\scalebox{0.5}{\includegraphics{minus.eps}}}}
%\caption{Subtraction of Two Sets}
\begin{picture}(160,100)(20,0)
\thinlines
\put(20,20){\makebox(50,25){$A - B$}}
\put(20,20){\line(0,1){25}}
\put(20,20){\line(1,0){50}}
\put(20,45){\line(1,0){50}}
\put(70,20){\line(0,1){5}}
\put(70,45){\line(0,-1){5}}
\put(10,10){\makebox{(b) Unexpected result}}

\put(140,20){\framebox(50,25){$A - B$}}
\put(130,10){\makebox{(c) Expected result}}

\put(70,70){\makebox(50,25){$A$}}
\put(120,75){\makebox(30,15){$B$}}
\put(70,70){\framebox(50,25){}}
\put(120,75){\line(1,0){30}}
\put(120,90){\line(1,0){30}}
\put(150,75){\line(0,1){15}}
\put(60,60){\makebox{(a) Original sets}}
\end{picture}
\caption{Subtraction of two sets}
\label{minus}
\end{figure}
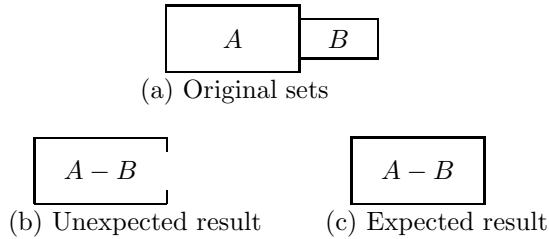

The theoretical role for $\mathcal{R}(\Re^3)$ was introduced into
geometric computing to correct the unexpected output seen from
combinatorial operations on geometric sets \cite{Requicha80}.  For
instance, consider the two dimensional illustration shown in
Figure~\ref{minus}.  The original operands of $A$ and $B$ are indicated in
Figure~\ref{minus}(a).  The unexpected output is shown in
Figure~\ref{minus}(b), where the expected result would have been what is
shown in Figure~\ref{minus}(c).

The phenomenon shown in Figure~\ref{minus}(b) was informally described as
``dangling edges'' \cite{VR77}.  The formalism that was proposed to
eliminate this behavior was that geometric combinatorial algorithms should
accept only regular closed sets as input and then execute the Boolean
operations of meet, join and complementation on these operands, thereby
creating only regular closed sets as output \cite{TilR80}.  The intent was
to eliminate ``dangling edges'' and, in principle, this should have been
sufficient\footnote{The subtraction operation between two sets, shown as
$A - B$ in Figure~\ref{minus}, is not {\em specifically} a Boolean
operation.  However, the use of $A - B$ should be understood to be
conveniently shortened notation equivalent to the operations $A \wedge
B^{\prime}$, where $B^{\prime}$ represents the standard Boolean operation
of complementation on the operand $B$.}.  However, each operand also has a
geometric representation that depends upon the approximation methods used
to compute the results. This additional subtlety raises issues in both
theory and computation.

For this short article, only a brief literature review will be presented.
An earlier survey on topology in computer-aided geometric design
\cite{PeRoD94} is recommended as introductory material for topologists.
The texts \cite{Ho89,PM-book} discuss the integration of computational
geometry, shape modeling and topology.  The subject of computational
topology is still a nascent and emerging sub-discipline.  This article
focuses upon the authors' particular perspective in its development. The
first use of the terminology ``computational topology'' appears to have
been in the dissertation of M.~M\"antyl\"a \cite{MM83}. Further
contemporary views can be gained from the following web sites
\cite{Top-manifesto,Edels-url,Z-url}.

Additional perspective can be gained by understanding the broader context
in which topology has already been successfully applied to computer
science.  We mention two particularly notable successes.  The first is the
use of non-Hausdorff topology by Kopperman, Meyer, Kong, Rosenfeld, Smyth
and Herman in digital topology for computer graphics and image processing.  
A good overview is readily available \cite{KoKoMe91}.  Similarly, the work
of Mislove, Reed and Roscoe on domains explores variants of limits for
fundamental algorithmic and programming language studies, continuing the
expressive power and the broad applicability of the language of topology
in denotational semantics and concurrent programming.  The monograph
\cite{Reed91} is recommended as an introduction.

\section{Theory versus computation}

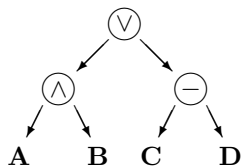
\begin{figure}[h]
%\centerline{\rotatebox{0}{\scalebox{0.5}{\includegraphics{tree.eps}}}}
%\setlength{GapDepth}{5mm}
%\setlength{GapWidth}{5mm}
%\setlength{GapWidth}{5mm}
\begin{picture}(110,65)
\put(50,50){\makebox(12,12){$\vee$}}
\put(56,56){\circle{12}}
\put(50,50){\vector(-1,-1){12}}
\put(62,50){\vector(1,-1){12}}
\put(25,25){\makebox(12,12){$\wedge$}}
\put(31,31){\circle{12}}
\put(25,25){\vector(-1,-2){6}}
\put(37,25){\vector(1,-2){6}}
\put(75,25){\makebox(12,12){$-$}}
\put(81,31){\circle{12}}
\put(75,25){\vector(-1,-2){6}}
\put(87,25){\vector(1,-2){6}}
\put(10,0){\makebox(12,12){$\mathbf{A}$}}
\put(40,0){\makebox(12,12){$\mathbf{B}$}}
\put(60,0){\makebox(12,12){$\mathbf{C}$}}
\put(90,0){\makebox(12,12){$\mathbf{D}$}}
\end{picture}
\caption{Tree for 
$(\mathbf{A} \wedge \mathbf{B}) \vee (\mathbf{C} - \mathbf{D})$.} 
\label{tree}
\end{figure}

One elegant computational representation for the combinatorial operators
is to assign each object a symbol and then to indicate operations in a
tree referencing those symbols.  For instance, such a tree structure could
be as depicted in Figure~\ref{tree}.

At this level of abstraction, the mathematical theory and the
computational representation are completely consistent, and this
representation became known as Constructive Solid Geometry (CSG).
Difficulties arose in instantiating the basic geometric information that
is represented by the operands at the leaf nodes and, sometimes, in
computing geometric representations at the internal nodes of the tree. In
CSG, the leaf nodes are restricted to a small set of specific geometric
objects, known as primitives. A typical collection of primitives might
consist of spheres, parallelepipeds, tori and right circular cylinders.  
The critical geometric algorithm underlying each Boolean operation is the
pairwise intersection between the operands.

As the boundary of each of these primitives can be represented by linear
or quadratic polynomials, the needed intersection between each pair of
primitives was relatively simple and numerically stable, for most cases
considered, although specific intersections could be problematical.  For
instance, suppose two cylinders of identical radius and height were
created and then positioned so that the bottom of one cylinder was
co-incident with the top of the other cylinder.  This special case was
specifically considered in most intersection algorithms and could usually
be processed without problem.  However, if one then rotated the top
cylinder a fraction of a degree about its axis (so that the planar
co-incidence remained intact) many software systems would fail to produce
any output for this problem, sometimes even causing a catastrophic program
failure.  This particular problem became a celebrated test case and most
systems developed {\it ad hoc} methods to solve this cylindrical
intersection problem.  Yet, this was just avoiding the more serious issue
of the fragile theoretical foundations for many intersection algorithms.  
People using CSG systems became sensitive to their limitations and
continued to use them effectively by avoiding these challenging
circumstances, although the work-arounds were often tedious to execute.

The imperative, largely initiated by the aerospace and automotive
industries, to model objects using polynomials of much higher degree than
quadratic created a movement away from CSG systems.  The alternative
format was to represent compact elements of $\mathcal{ R}(\Re^3)$ by their
boundaries, and this became known as the ``boundary representation''
approach, or ``B-rep'' for short.  This has become the dominant mode
today.  Again, within this clean conceptual overview, the realities of
computation pose some subtle problems.  In most industrial practice, the
modeling paradigm was further restricted so that the boundary of an object
was a 2-manifold without boundary.  However, it was difficult to create
computer modeling tools that could globally define 2-manifolds without
boundary, though there existed excellent tools for creating subsets of
these 2-manifolds. For example, computational tools for creating splines
were becoming prevalent.  Again, in principle, if each such spline subset
was created with its boundaries, then the subsets could be joined along
shared boundary elements to form a topological complex \cite{HY} for the
bounding 2-manifold without boundary.

The inherent computational difficulty was to separately create two spline
patches, each being a manifold with boundary, so that the corresponding
boundary curves were identical and could be exactly shared between the
patches. In some situations, algorithms for fitting spline patches were
used successfully.  In other cases, patches have been slightly enlarged
and intersected so as to obtain improved fits. Indeed, such intersections
are well-defined in pure mathematics, but, again, approximation in
computation poses subtle variations from that theory, as described in the
next section on pairwise surface intersection.

\section{Subtleties of pairwise spline surface intersection}

It is well-known that unwanted gaps between spline surfaces or
self-intersections within intended manifolds often appear as unwanted
artifacts of various implemented intersection algorithms
\cite{Fa-siam-news}.  The mismatch between approximate geometry and exact
topology has historically caused reliability problems in graphics, CAD,
and engineering analyses---drawing the attention of both academia and
industry. The severity of the problem increases with the complexity of the
geometric data represented, both from high-degree nonlinearity and from
the intricate interdependence of shape elements that should, but do not,
fit together according to the specified topological adjacency information.

\begin{figure}[ht]
%\centerline{\scalebox{.8}{\includegraphics{int-F-G-p77.eps}}}
\centerline{\scalebox{.5}{\includegraphics{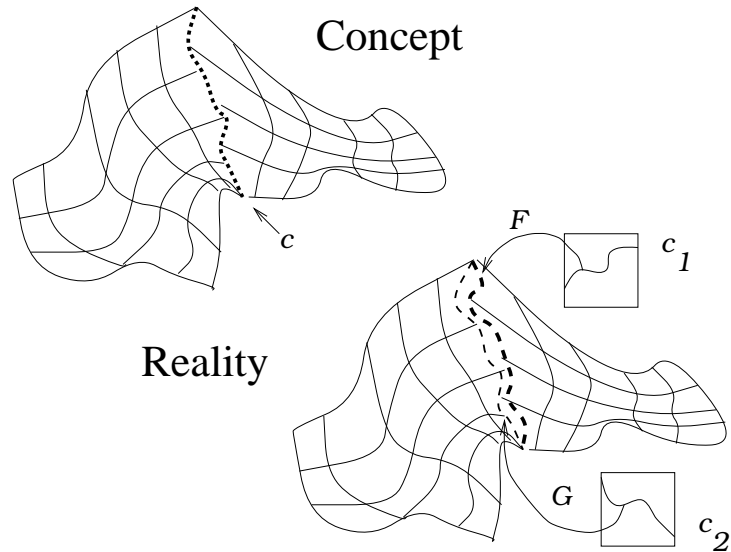}}}
\caption{Joining operations for geometric objects}
\label{gaps}
\end{figure}

The conceptual view of these joining operations is illustrated in
Figure~\ref{gaps}, with an intersection curve\footnote{We focus on the
generic case of an intersection curve, although isolated points and
co-incident areas can also arise, with similar complications.} denoted as
$c$. But this picture of $c$ only exemplifies the idealized, exact
intersection curve.  For practical computations, an approximation of $c$
is often created \cite{G-siam-review} and, in many systems, an
intersection curve will be approximated twice; once within the parametric
domain of one of the intersecting surfaces and then again within the
parametric domain of the other.  These approximations are labeled as $c_1$
and $c_2$ in Figure~\ref{gaps}.  The spline functions from $[0,1]^2$ into
$\Re^3$ then also rely upon algorithmic evaluation of these approximated
intersection curves, as indicated by $F$ and $G$ in Figure~\ref{gaps}.  
It is virtually certain that those evaluations will not be exactly equal
in $\Re^3$.

The mismatch between concept and reality depicted in Figure~\ref{gaps}
creates ambiguity, as the intersection representation is sometimes
considered as a unique set, from the symbolic topological view, and at
other times as two approximating sets, from the geometric view.  

\section{Specific progress}

To resolve this ambiguity, we are investigating richer representations.  
We introduce a neighborhood of the true, but unknown, intersection set.
This neighborhood is created from newly determined rigorous upper bounds
on the error incurred during efficient intersection approximations.  To
date, it has been convenient to create these neighborhoods as tubular
neighborhoods \cite{Hi76}, but broader generalizations seem to be
possible.

\subsection*{The role of topological equivalence}

\begin{figure}[htbp]
%\centerline{\scalebox{.5}{\includegraphics{2-knots-70p-40font.eps}}}
\centerline{\scalebox{.3333}{\includegraphics{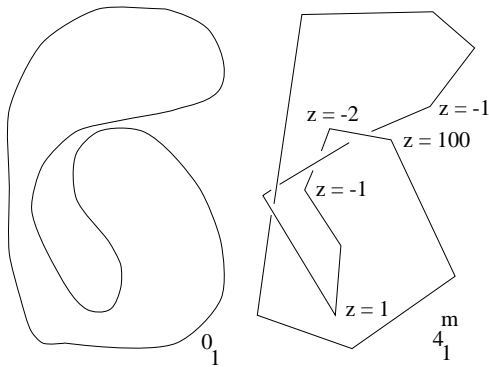}}}
\caption{Nonequivalent knots}
\label{knot}
\end{figure}

Before introducing these bounding neighborhoods, we discuss the meaning of
topological equivalence between an object and its approximations.  We have
proposed the use of ambient isotopy for this topological equivalence
versus the more traditional equivalence by homeomorphism, as is explained
more fully in our publications
\cite{APR03,ADPS95,APS00,SP03,SPB-CAD,St93}.  Intuitively, two closed
curves will not be ambient isotopic if they form different knots.  
Figure~\ref{knot} shows two simple homeomorphic space curves, where the
piecewise linear (PL) curve is an approximation of the smooth curve.  
However, these curves are not ambient isotopic, because they depict
different knots\footnote{The different knot classifications of $0_1$ and
$4^m_1$ are indicated in Figure~\ref{knot}.}, with the smooth curve
illustrating the simplest knot, known as the unknot.  In the right half of
Figure~\ref{knot} the $z$ coordinates of some vertices are specifically
indicated to emphasize the knot crossings in $\Re^3$ (All other end points
have $z = 0$).  All end points of the line segments in the PL
approximation are also points on the original curve. Having this knotted
curve as an approximant to the original unknot would be undesirable as
output from a curve approximation algorithm, particularly for applications
in graphics and engineering simulations. Similar pathologies can happen in
approximating 2-manifolds, both with and without boundary, but results
\cite{APR03,APS00,NMP-TS-98,SP03,SPB-CAD} summarized here can prevent
these difficulties by appropriately constraining the approximations
produced.

In response to the example of Figure~\ref{knot}, a theorem was published
that provided for ambient isotopic PL approximations of 1-manifolds
\cite{NMP-TS-98}.  The proof utilizes ``pipe surfaces'' from classical
differential geometry \cite{Monge} to build an appropriately small tubular
neighborhood such that if the PL approximant is constrained to lie within
the constructed neighborhood, then the PL approximant is ambient isotopic
to the original curve.  While the techniques of tubular neighborhoods are
well known in differential topology, the relevant theorems usually state
only the existence of an ambient isotopy. Additional work in applied
mathematics was needed to complement the theorems with constructive
formulations from which to obtain effective procedures and algorithms.  
Such extensions were successfully accomplished for both non-singular
compact, orientable 1-manifolds and 2-manifolds (with or without boundary)
\cite{Abeet04,APR03,NMP-TS-98,SP03,SPB-CAD} by I-TANGO team members (and
their co-authors).  We used geometric characteristics to compute a {\em
specific} upper bound on the size of a tubular neighborhood and to then
specify a {\em particular} isotopy.

The importance of this topological equivalence class extends beyond the
manifolds described to the geometric models created as compact subsets of
$\mathcal{R} (\Re^3)$.  Consider that a well-defined topological complex
could be created from 2-manifolds with boundary if the difficulties along
the intersection boundaries could be solved.  However, the non-uniqueness
of the geometric representation of the intersection sets seems to pose an
intractable problem to creation of a single, well-defined topological
complex.  Our alternative approach is to find a neighborhood within which
it can be proven that the true intersection curve lies.  Since the
construction of these bounding neighborhoods is dependant upon the
specific intersection algorithm used, some further details of the two
intersection algorithms used within the I-TANGO project are presented.  
While the neighborhood construction details will vary with the specific
intersection algorithm chosen, the following two neighborhood
constructions were chosen both for their carefully defined error bounds
and for their potential for generalization.

\subsection*{Error bounds for topology from Taylor's theorem}

First, we present the Grandine-Klein (GK) intersection algorithm
\cite{GK97}. Referring to Figure~\ref{gaps}, we note that the GK algorithm
bases its error bounds on well-established numerical techniques in
differential algebraic equations (DAE).  While these DAE techniques
provide rigorous error bounds, these bounds are expressed within the
parameter space $[0,1]^2$, which serves as the domain of the spline
functions (indicated as $F$ and $G$ in Figure~\ref{gaps}).  The code
implementing the GK algorithm then has an interface that allows the user
to specify an upper bound $\epsilon$ for this error in parameter space and
the algorithm provides guarantees for meeting this error bound.  However,
the typical end user is often unaware of the role of this parametric
domain, so selection of this parametric space error bound has often relied
upon heuristics.  It would be more convenient for the user to be able to
specify an error bound within $\Re^3$.  One accomplishment within the
I-TANGO project has been to demonstrate a mathematical relation
\cite{MPS03} between the error bounds in $\Re^3$ and $[0,1]^2$, following
from a straightforward application of Taylor's Theorem in two dimensions
\cite[p.~200]{Buck}.  The conversion between these error bounds has been
implemented in a pre-processing interface to the GK algorithm and this new
interface has been tested to be reliable, efficient and user-friendly.

Using the notation from Figure~\ref{gaps} for the spline function $F$, 
Taylor's Theorem provides a bound on the error of $F$ evaluated at a
particular point $(u,v)$ versus $F$ evaluated at a point $(u_0, v_0)$,
where $(u,v)$ and $(u_0, v_0)$ are within a sufficiently small neighborhood. 
This sufficiently small neighborhood will have radius given by the value
in the parametric domain $[0,1]^2$ which was denoted as $\epsilon$
in the previous paragraph.  Then it follows \cite{MPS03}, with 
$\| \cdot \|$ being any convenient vector norm, that
\begin{equation*}
\|F(u,v) - F(u_0,v_0) \| \leq \epsilon M
\end{equation*}
for any $M$ satisfying
\begin{equation*}
\left\| \frac{\partial F}{\partial u}(u^{*},v^{*}) \right\| 
+ \left\| \frac{\partial F}{\partial v}(u^{*},v^{*}) \right\| 
\leq M,
\end{equation*}
for some point $[u^*,v^*]$ on the line segment joining
$[u,v]$ and $[u_1,v_1]$.

For the single spline $F$, let $\gamma(F)$ be an upper bound for the 
acceptable error in $\Re^3$ between the true intersection curve $c$ and 
one of its approximants $F(c_1)$.  In order to guarantee that this error 
is sufficiently small, it is sufficient
$$\epsilon M \leq \gamma(F),$$
where an upper bound for $M$ can be computed by using any standard 
technique for obtaining the maximums of the partials 
$\frac{\partial F}{\partial u}$ and $\frac{\partial F}{\partial v}$.  
For $G$, a similar relation between $\gamma(G)$ and $\epsilon$ 
exists\footnote{This error bound assumed that the error due to 
algorithmic truncation within the numerical DAE methods dominated any 
other computational errors.}.  

Then it is clear that a neighborhood can be defined that contains
the true intersection curve $c$ and both of its approximants.
Let $F(c_1)$ denote the image of $c_1$ under $F$ and similarly, let
$G(c_2)$ denote the image of $c_2$ under $G$.  
Let $N_{\gamma(F)}(F(c_1))$ be a tubular neighborhood of radius 
$\gamma(F)$ about $F(c_1)$, where $c_1$ has been generated from the GK 
intersector to satisfy the inequality presented in the previous 
paragraph.  Similarly, define $N_{\gamma(G)}(G(c_2))$.  Then, let 
\[
N(c) = N_{\gamma(F)}(F(c_1)) \cup N_{\gamma(G)}(G(c_2)).
\]

It is clear that $N(c)$ is a neighborhood of $c$, which contains both
of its approximants, $F(c_1)$ and $G(c_2)$.  However, there is both a
theoretical and computational limitation to this approach.  

\begin{itemize}
\item There is no theoretical guarantee that either approximant is
topologically equivalent to $c$, and
\item Any practical computation of $N(c)$ would depend upon an
accurate computation of the set 
\mbox{$N_{\gamma(F)}(F(c_1)) \cap N_{\gamma(G)}(G(c_2))$},
which is likely to be as difficult as the original computation of $F \cap G$.
\end{itemize}

While the above bounds are often quite acceptable in practice to compute a
reasonable approximant, further research has been completed into alternate
methods to give guarantees of topological equivalence within a
computationally acceptable neighborhood of the intersection set, as
reported in the next subsection.

\subsection*{Integrating error bounds and topology via interval solids}

Recent work by Sakkalis, Shen and Patrikalakis \cite{SSP01} emphasized that
the numeric input to any intersection algorithm has an initial
approximation in the co-ordinates used to represent points in $\Re^3$,
leading to their use of interval arithmetic \cite{PM-book}.  The basic
idea behind interval arithmetic is that any operation on a real value $v$
is replaced by an operation of an interval of the form $[a,b]$, where $a,
b \in \Re$ and $a < v < b$.  The result of any such interval operation is
an interval, which is guaranteed to contain the true result of the
operation on $v$.  This led naturally to the concept of an {\em interval
solid} and some of its fundamental topological and geometric properties
were then proven, as summarized below.

Throughout this section, a {\em box} is a rectangular, closed 
parallelepiped in $\Re^3$ with positive volume, whose edges are parallel
to the co-ordinate axes\footnote{Enclosures other than boxes are quite
possible and this is a subject of active research.}.  Let $F$ be a
non-empty, compact, connected 2-manifold without boundary.  Then the
Jordan Surface Separation Theorem asserts that the complement of $F$ in
$\Re^3$ has precisely two connected components, $F_I$, $F_O$; we may 
assume that $F_I$ is bounded and $F_O$ is unbounded.  
Let also $\mathcal{B} = \{ b_j, j \in J\}$ be a finite collection of boxes 
that satisfies the following conditions:
\begin{description}
\item[C1] 
$\{ \Int(b_j), j \in J \}$ is a cover of $F$.
\item[C2] 
Each member $b$ of $\mathcal{B}$ intersects $F$ generically; that is, 
$b \cap F$ is a non-empty closed disk that separates $b$ into two 
(closed) balls, $B_b^+$ and $B_b^-$, with $B_b^+$, ($B_b^-$) lying in 
$F_I \cup F$ ($F_O \cup F$), respectively.
\item[C3] 
For any $b_{i}, b_{j} \in \mathcal{B}$, let $b_{ij}= b_i \cap b_j$. 
If $\Int(b_i) \cap \Int(b_j) \neq \emptyset$, then $b_{ij}$ is also a box 
which satisfies {\bf C2}.
\end{description}

\begin{figure}[htbp]
%\centerline{\scalebox{.5}{\includegraphics{c2-c3.eps}}}
\centerline{\scalebox{.3333}{\includegraphics{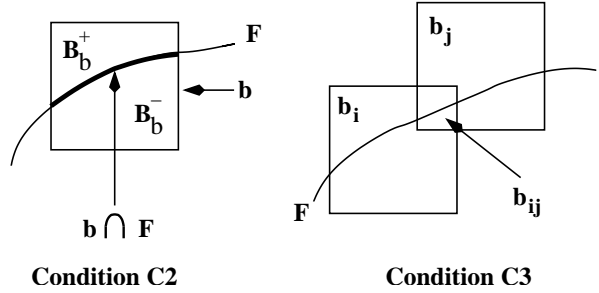}}}
\caption{2D versions of conditions {\bf C2} and {\bf C3}}
\label{cond4-5}
\end{figure}

Notice that condition {\bf C2} indicates that every $b \in \mathcal{B}$
intersects $F$ in a natural way (see Figure~\ref{cond4-5}).

%\begin{picture}(100,100)
%\thinlines
%\put(12,36){\makebox(12,12){$b$}}
%\put(62,11){\makebox(12,12){$B_b^-$}}
%\put(62,36){\makebox(12,12){$B_b^+$}}
%\put(10,10){\framebox(40,40){}}
%\put(-10,10){\makebox{$F$}}
%\thicklines
%\put(0,0){\qbezier(10,20)(22,33)(50,34)}
%\thinlines
%\put(0,0){\qbezier(0,10)(20,35)(60,35)}
%\end{picture}
%\begin{picture}(100,100)
%\put(10,10){\framebox(40,40){}}
%\put(35,35){\framebox(40,40){}}
%\put(0,0){\qbezier(0,15)(30,50)(85,55)}
%\end{picture}

The following result summarizes several previously appearing results,
where a solid is defined to be a non-empty compact, regular closed subset
of $\Re^3$.

\begin{thm}[{\cite[Corollary 2.1, p.~165]{SSP01}}]
\label{intsolid} 
If $F$ is connected and $\mathcal{B}$ satisfies {\bf C1}--{\bf C3}, 
then $F \cup \bigcup_{j \in J} b_j$ is a solid.
\end{thm}

Bisceglio, Peters and Sakkalis \cite{SP03,SPB-CAD} have recently given 
sufficient conditions to show when the boundary of an interval solid is ambient
isotopic to the well-formed solid that it is approximating, as described 
in the following theorem.
For a positive number $\delta$, define the open set
\[ F(\delta)=\{x \in \Re^3 \,|\, \hbox{$D(x, F) < \delta$} \}, \]
where $D(x, F) = \inf \{d(x,y) \,|\, y \in F\}$, with $d$ being the usual
Euclidean metric in $\Re^3$.

\begin{figure}[htbp]
%\centerline{\scalebox{.9}{\includegraphics{proper.eps}}}
\centerline{\scalebox{.4}{\includegraphics{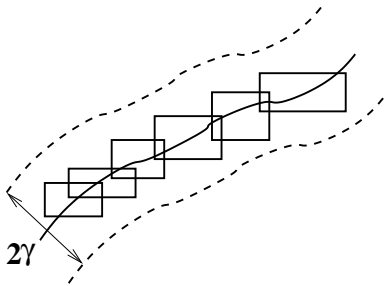}}}
\caption{2D version of proper subset condition }
\label{proper-sub}
\end{figure}

\begin{thm}
\label{ai-sp}
Let $F$ be a connected 2-manifold without boundary. 
For each $\epsilon > 0$, there exists $\delta$, with $0 < \delta < \rho$ 
so that whenever a family of boxes $\mathcal{B}$ satisfies conditions 
{\bf C1}--{\bf C3}, and for each $b$ of $\mathcal{B}$, $b$ is a  subset 
of $F(\delta)$ \emph{(Please see Figure~\ref{proper-sub})} then, for 
$S=F \cup F_I$ and $S^\mathcal{B} = S \cup \bigcup_{j \in J} b_j$, the 
sets $F$ and $\partial S^\mathcal{B}$ are $\epsilon$-isotopic with 
compact support.  
Hence, they are also ambient isotopic.
\end{thm}

The quoted theorem depends upon results from Bing's book on PL topology
\cite[p.~214]{Bing83}, and related literature \cite{Ki59}, as is explained
in full \cite{SP03,SPB-CAD}. The proof shows that normals to $F$ do not
intersect within the constructed tubular neighborhood, as is illustrated
by the depiction of its planar cross-section in Figure~\ref{proper-sub}.

If the boxes containing the true intersection curve can be made
sufficiently small so that each such box fits inside $F(\rho)$, then the
resultant intersection neighborhood will contain an object that is both
close to the true solid and is ambient isotopic to it. Considerable
success in meeting these constraints has already been achieved
\cite{PM-book} when two splines intersect transversally, while very recent
progress from the I-TANGO team for more subtle spline intersection pairs
is now under review \cite{MKMSP03}.  As a further note on integration,
results on root computations from interval arithmetic \cite{PM-book} are
used to provide estimates of initial starting points for the GK algorithm.

\subsection*{Work in Progress}

The previous discussion emphasizes results from the I-TANGO project that
have already appeared in the literature, whereas the remarks in this
section are intended to provide some indication of related results that
are expected to appear soon.

Recall that the result of the GK algorithm will be a set of spline patches
that do not fit together precisely along the approximation of their boundary. 
Recent work by Andersson, Stewart and Zidani \cite{ASZ03} uses the Whitney
Extension Theorem to formulate a conceptual model of how imperfectly
fitting patches might be perturbed to form a topological complex which is
the boundary of a non-empty compact element of $\mathcal{R}(\Re^3)$.
Furthermore, the proposed process extends the patches under Lipschitz
mappings so that the resulting element of $\mathcal{R}(\Re^3)$ can be
shown to lie within rigorous error bounds of the given geometric input
data. Our initial results on ambient isotopic approximations of
\mbox{2-manifolds} have followed the prevailing simplification of
considering only those manifolds without boundaries
\cite{APR03,SP03,SPB-CAD}, but our recent results have been for the more
technically challenging cases of 2-manifolds with boundaries
\cite{Abeet04}.  Our team is also completing work on spline intersections
with multiple roots \cite{MKMSP03}, to improve the rather loose error
bounds known to-date.

Interestingly, the interval solids were not initially intended as a means
to integrate issues of topological equivalence and approximations in
intersection algorithms.  Yet it is a hallmark of our project that these
concepts are merging, as discussed above.  Note, further, that the
intervals include both approximation errors from truncation of numerical
processes within intersection algorithms as well as the approximations
that arise by using a finite set of floating point numbers in computation
as an approximation of the reals.  How floating point approximations
effect error bounds for intersection computations is another important
issue this project considers.  Attention has focused upon the impact of
floating point arithmetic on polynomial computations near multiple roots
\cite{HPSS03}.  This is an important issue because the intersection
algorithms have to search for roots of the system of polynomial equations
representing the intersecting surfaces and exact arithmetic computations
have not been shown to be practical.

\section{Conclusions and future work}

Topology and computer science are finding common interest in the emerging
area of computational topology.  Various branches of pure topology
(point-set topology, differential topology, low-dimensional topology,
$\ldots$) can make important contributions to establishing the appropriate
theoretical foundations.  A fundamentally new perspective arises from the
role that computational approximations should play in the reformulation of
central topological concepts.  The I-TANGO project is a research effort
concentrating upon these issues specifically with respect to surface
intersections within $\Re^3$.  Considerable progress has already been
made, but many questions also remain open, as is summarized in this
article.

\end{document}